\documentclass[12pt]{article}
\catcode`\@=11
\newtheorem{theorem}{Theorem}
\newtheorem{lemma}{Lemma}
\newtheorem{proposition}{Proposition}

\newtheorem{example}{Example}

\newtheorem{corollary}{Corollary}
\def\demo{\noindent{\bf Proof .-}}
\def\section{\@startsection {section}{1}{\z@}{-3.5ex plus -1ex
minus-.2ex}{2.3ex plus .2ex}{\normalsize\bf}}
\def\codim{{\rm codim}\,}

\def\bn{\hbox{\it I\hskip -2pt N}}
\def\bz{\hbox{\it Z\hskip -4pt Z}}

\newcommand{\het}{H_{\rm et}}
\newcommand{\hc}{H_{\rm c}}
\begin{document}
\begin{center}
{\Large\bf \textsc{On toric varieties which are almost set-theoretic complete intersections}}
\end{center}
\vskip.5truecm
\begin{center}
{Margherita Barile\footnote{Partially supported by the Italian Ministry of Education, University and Research.}\\ Dipartimento di Matematica, Universit\`{a} di Bari,Via E. Orabona 4,\\70125 Bari, Italy}
\end{center}
\vskip1truecm
\noindent
{\bf Abstract} We describe a class of affine toric varieties $V$ that are set-theoretical\-ly minimally defined by $\codim V+1$ binomial equations over fields of any characteristic. \vskip0.5truecm
\noindent

\section*{Introduction}
Let $K$ be an algebraically closed field. Let $V$ be an affine variety, $V\subset K^N$, where $N$ is the smallest possible. The arithmetical rank (ara) of $V$ is defined as the least number of equations that are needed to define $V$ set-theoretically as a subvariety of $K^N$. In general we have that ara\,$V\geq\codim\,V$. If equality holds, $V$ is called a {\it set-theoretic complete intersection}. If ara\,$V\leq\codim\,V+1$, $V$ is called an {\it almost set-theoretic complete intersection}. The problem of determining when a toric variety is (almost) set-theoretic complete intersection is open in general; it includes the still unsettled case of monomial curves, which was intensively studied over the last three decades (see, e.g., \cite{Br}, \cite{El}, \cite{E}, \cite{K}, \cite{Mo}, \cite{RV}, \cite{T1}). In \cite{BMT1} it was shown that, whenever $V$ is a simplicial toric variety with a full parametrization, then ara\,$V=\codim\,V$ in all positive characteristics,  whereas ara\,$V\leq\codim\,V+1$ in characteristic zero. A different class of toric varieties with the same property has been recently described in \cite{B1}. In this paper we present a class of toric varieties $V$ of arbitrary codimension greater than 2, such that ara\,$V=\codim\,V+1$ in all characteristics: each variety $V$ is defined, over any field, by the vanishing of the same set of $\codim\,V+1$ binomials, even though these do not generate the defining ideal $I(V)\subset K[x_1,\dots, x_N]$. It seems that no example of this kind has ever appeared in the recent literature.

\section{Preliminaries} Let $K$ be an algebraically closed field.
 Let $n\geq3$ be an integer and let ${\bf e}_1,\dots,{\bf e}_n$ be the standard basis of $\bz^n$. Set $N=2n$ and consider the following subset of $\bn^{N}$:
\begin{eqnarray*}T&=&
\{d_1{\bf e}_1,\dots, d_{n-1}{\bf e}_{n-1}, {\bf e}_n,\\
 &&f_1{\bf e}_1+g_1{\bf e}_n,\dots, f_{n-1}{\bf e}_{n-1}+g_{n-1}{\bf e}_n, h_1{\bf e}_1+\cdots+h_{n-1}{\bf e}_{n-1}\},
\end{eqnarray*}
\noindent 
where $d_1,\dots, d_{n-1}, f_1,\dots, f_{n-1}, g_1,\dots, g_{n-1}, h_1,\dots, h_{n-1}$ are all positive integers such that 
\begin{equation}\label{1} \gcd(d_i,f_i)=1\mbox{ for all }i=1,\dots, n-1,\end{equation}
\begin{equation}\label{1b} \gcd(d_i,h_i)=1\mbox{ for all }i=1,\dots, n-1,\end{equation}
\noindent and
\begin{equation}\label{2} \gcd(d_i,d_j)=1\mbox{ for all }i,j=1,\dots, n-1, i\neq j.\end{equation}
\noindent
Suppose, moreover, that there are two distinct (positive) primes $p$ and $q$ such that $p$ divides $d_i$ and $q$ divides $d_j$ for some indices $i$ and $j$.  With $T$ we can associate the affine variety $V\subset K^N$ admitting the following parametrization
$$V:\left\{
\begin{array}{rcl}
x_1&=&u_1^{d_1}\\
&\vdots&\\
x_{n-1}&=&u_n^{d_{n-1}}\\
x_n&=&u_n\\
y_1&=&u_1^{f_1}u_n^{g_1}\\
&\vdots&\\
y_{n-1}&=&u_{n-1}^{f_{n-1}}u_n^{g_{n-1}}\\
y_n&=&u_1^{h_1}\cdots u_{n-1}^{h_{n-1}}\\
\end{array}\right.$$
\noindent
We have that $\codim V=n$.
The polynomials in the defining ideal $I(V)$ of $V$  are the linear combinations of binomials
$$B_{\alpha_1,\dots,\alpha_n,\beta_1,\dots,\beta_n}=
x_1^{\alpha_1^+}\cdots x_n^{\alpha_n^+}y_1^{\beta_1^+}\cdots y_n^{\beta_n^+}-
x_1^{\alpha_1^-}\cdots x_n^{\alpha_n^-}y_1^{\beta_1^-}\cdots y_n^{\beta_n^-}$$
where $({\alpha_1,\dots,\alpha_n,\beta_1,\dots,\beta_n})\in\bz^N\setminus\{\bf{0}\}$  is such that
\begin{eqnarray*}&&\alpha_1d_1{\bf e}_1+\cdots+\alpha_{n-1}d_{n-1}{\bf e}_{n-1}+\alpha_n{\bf e}_n+\\
&&\beta_1(f_1{\bf e}_1+g_1{\bf e}_n)+\cdots+\beta_{n-1}(f_{n-1}{\bf e}_{n-1}+g_{n-1}{\bf e}_n)+\\
&&\beta_n(h_1{\bf e}_1+\cdots+h_{n-1}{\bf e}_{n-1})=\bf{0}\qquad\qquad\qquad\quad\qquad\qquad\qquad\qquad\qquad\qquad(\ast)
\end{eqnarray*}
\noindent
and $\alpha_i^+=\max\{\alpha_i, 0\}$, $\alpha_i^-=\max\{-\alpha_i, 0\}$, $\beta_i^+=\max\{\beta_i, 0\}$ and  $\beta_i^-=\max\{-\beta_i, 0\}$.
There is a one-to-one correspondence between the set of binomials in $I(V)$ and the set of relations $(\ast)$ between the elements of
$T$.  Our aim is to prove  the following
\begin{theorem}\label{main} ara\,$V=n+1$.
\end{theorem}
This will be done by proving the two inequalities separately, in Sections 2 and 4 respectively.
\section{The defining equations}
In this section we explicitly exhibit $n+1$ binomial equations which define $V$ set-theoretically over any field $K$. \par\noindent
From (\ref{1}) it follows that there are $\alpha_1,\dots, \alpha_{n-1},\beta_1\,\dots,\beta_{n-1}\in\bz$ such that
\begin{equation}\label{bezout}h_i=\alpha_id_i+\beta_if_i\qquad\mbox{ for all }i=1,\dots, n-1.\end{equation}
\noindent
Then 
$$\sum_{i=1}^{n-1}\alpha_id_i{\bf e}_i-\left(\sum_{i=1}^{n-1}\beta_ig_i\right){\bf e}_n+\sum_{i=1}^{n-1}\beta_i(f_i{\bf e}_i+g_i{\bf e}_n)-\sum_{i=1}^{n-1}h_i{\bf e}_i=0.$$
\noindent
Hence the binomial
$$G=B_{(\alpha_1,\dots, \alpha_{n-1}, -\sum_{i=1}^{n-1}\beta_ig_i, \beta_1,\dots,\beta_{n-1}, -1)}$$
belongs to $I(V)$. Let
$$F=y_n^{d_1\cdots d_{n-1}}-x_1^{h_1d_2\cdots d_{n-1}}\cdots x_i^{d_1\cdots d_{i-1}h_{i}d_{i+1}\cdots d_n}\cdots x_{n-1}^{d_1\cdots d_{n-2}h_{n-1}},$$
\noindent
and
$$F_i=y_i^{d_i}-x_i^{f_i}x_n^{d_ig_i}\qquad\mbox{ for all }i=1,\dots, n-1.$$
\noindent
An easy computation yields that $F, F_1,\dots, F_{n-1}\in I(V)$.
\begin{proposition}\label{equations} Variety $V$ is set-theoretically defined by
$$F_1=\cdots=F_{n-1}=F=G=0.$$
\end{proposition}
\demo We only have to prove that every ${\bf w}\in K^N$ fulfilling the given equations belongs to $V$. So let ${\bf w}=(\bar x_1,\dots, \bar x_n,\bar y_1,\dots,\bar y_n)\in K^N$ be such that 
\begin{eqnarray}F_1({\bf w})=\cdots= F_{n-1}({\bf w})&=&0\label{F1}\\
F({\bf w})&=&0\label{F}\\
G({\bf w})&=&0\label{G}
\end{eqnarray}
\noindent
Let $u_1,\dots, u_n\in K$ be such that $\bar x_i=u_i^{d_i}$ for all $i=1,\dots, n-1$, and set $u_n=\bar x_n$. We show that for a suitable choice of $u_1,\dots, u_{n-1}$,   we have 
\begin{equation}\label{a}\bar y_i=u_i^{f_i}u_n^{g_i}\quad\mbox{ for all }i=1,\dots, n-1,\end{equation}
\noindent and
\begin{equation}\label{b}\bar y_n=u_1^{h_1}\cdots u_{n-1}^{h_{n-1}}.\end{equation}
\noindent
We first show that, up to replacing some $u_i$ with another $d_i$-th root of $\bar x_i$,  (\ref{a}) is fulfilled.   From (\ref{F}) we deduce that $\bar y_n^{d_1\cdots d_{n-1}}=u_1^{h_1d_1\cdots d_{n-1}}\cdots u_{n-1}^{h_{n-1}d_1\cdots d_{n-1}}$, so that
$$\qquad\qquad\qquad\qquad\qquad\bar y_n=u_1^{h_1}\cdots u_{n-1}^{h_{n-1}}\eta,\qquad\qquad\qquad\qquad\qquad\qquad\quad(9)'$$
\noindent
for some $\eta\in K$ such that $\eta^{d_1\cdots d_{n-1}}=1.$  For all $i=1,\dots, n-1$ let $\eta_i\in K$ be a primitive $d_i$-th root of unity. Then (\ref{2}) implies that $\eta'=\eta_1\cdots\eta_{n-1}$ is a primitive $d_1\cdots d_{n-1}$-th root of unity. Hence, for a suitable positive integer $s$, we have $\eta'^s=\eta$. Choose $u_i\eta_i^s$ instead of $u_i$, for all $i=1,\dots, n-1$. Then (9)$\prime$ turns into (\ref{b}). Condition (\ref{a}) is also fulfilled if $\bar x_n=0$. So suppose that $\bar x_n\ne0$. We prove that, up to a modification of $u_1,\dots, u_{n-1}$ which preserves (\ref{b}), condition (\ref{a}) is satisfied. From (\ref{F1}) we derive that, for all $i=1,\dots, n-1$, $\bar y_i^{d_i}=u_i^{d_if_i}u_n^{d_ig_i}$, so that  
$$\qquad\qquad\qquad\qquad\qquad\bar y_i=u_i^{f_i}u_n^{g_i}\theta_i,\qquad\qquad\qquad\qquad\qquad
\qquad\qquad\qquad(8)'$$
\noindent where $\theta_i\in K$ is some $d_i$-th root of unity. From (\ref{1}) it follows that $\theta_i=\omega_i^{f_i}$ where $\omega_i\in K$ is some $d_i$-th root of unity.  We replace $u_i$ by $u_i\omega_i$ for all $i=1,\dots, n-1$. Then (8)$\prime$ turns in to (\ref{a}). Finally, from (\ref{G}) it follows that $\bar y_n\bar x_n^{\beta_1g_1+\cdots+\beta_{n-1}g_{n-1}}=\bar x_1^{\alpha_1}\cdots \bar x_{n-1}^{\alpha_{n-1}}\bar y_1^{\beta_1}\cdots\bar y_{n-1}^{\beta_{n-1}}$, whence
$$\bar y_nu_n^{\beta_1g_1+\cdots+\beta_{n-1}g_{n-1}}=u_1^{\alpha_1d_1+\beta_1f_1}\cdots u_{n-1}^{\alpha_{n-1}d_{n-1}+\beta_{n-1}f_{n-1}}u_n^{\beta_1g_1+\cdots\beta_{n-1}g_{n-1}},$$
\noindent
which, being $u_n\ne0$, and in view of (\ref{bezout}), implies (\ref{b}).
This completes the proof. We have just proven that: 
\begin{corollary} ara\,$V\leq n+1$.
\end{corollary}
In general, even if $V$ is set-theoretically defined by $n+1$ binomial equations, $n+1$ binomials are not sufficient to generate
the defining ideal of $I(V)$. This is shown in the next
\begin{example}{\rm  Let $n=3$, and take $d_1=2$, $d_2=3$, $f_1=3$, $f_2=5$, $g_1=g_2=1$, $h_1=3$, $h_2=5$. The corresponding toric variety is
$$V:\left\{\begin{array}{rcl}
 x_1&=&u_1^2\\
x_2&=&u_2^3\\
x_3&=&u_3\\
y_1&=&u_1^3u_3\\
y_2&=&u_2^5u_3\\
y_3&=&u_1^3u_2^5\\
\end{array}\right.$$
A computation by CoCoA \cite{CoCoA} yields that the defining ideal $I(V)$ is minimally generated by the following 8 binomials:
$$y_1^2-x_1^3x_3^2,\qquad y_2^3-x_2^5x_3^3,\qquad y_3^6-x_1^9x_2^{10},\qquad x_3^2y_3-y_1y_2,$$
$$y_1y_3-x_1^3y_2,\qquad y_2y_3^2-x_1^3x_2^5x_3,\qquad y_2^2y_3-x_2^5x_3y_1,\qquad x_3y_3^3-x_1^3x_2^5y_1.$$
\noindent
Since $f_i=h_i$ for $i=1,2$,  in (\ref{bezout}) we can take $\alpha_i=0$ and $\beta_i=1$ for $i=1,2$. Then the first four binomials in the above list are $F_1, F_2, F$ and $G$ respectively. According to Proposition \ref{equations} they suffice to define $V$ set-theoretically. }
\end{example}

\section{Some cohomological results}
In the next section we shall complete the proof of Theorem \ref{main} using cohomological tools. In this section we provide the necessary preliminary lemmas on \'etale cohomology. We refer to \cite{M} for the basic notions. Let $K^{\ast}=K\setminus\{0\}$. 
\begin{lemma}\label{lemma1} Let $n$ be a positive integer, and let $r$ be an integer prime to char\,$K$.  Let $d_1,\dots, d_n$ be positive integers such that, for some index $i$,  $d_i$ and $r$ are not coprime, and consider the morphism of schemes
$$\delta:(K^{\ast})^n\rightarrow  (K^{\ast})^n$$
$$(u_1,\dots, u_{n})\mapsto (u_1^{d_1},\dots, u_{n}^{d_{n}})$$
\noindent 
Then the map 
$$\kappa_n:\hc^{n+1}((K^{\ast})^n)\to\hc^{n+1}((K^{\ast})^n)$$
\noindent
induced in cohomology with compact support is not injective.
\end{lemma}
\demo In the sequel $\het$ and $\hc$ will denote \'etale cohomology and cohomology with compact support with respect to $\bz/r\bz$; for the sake of simplicity we shall omit the coefficient group in this and in the next proofs. We proceed by induction on $n\geq1$. For $n=1$ we have the morphism 
$$\delta:K^{\ast}\rightarrow K^{\ast}$$
$$u_1\mapsto u_1^{d_1},$$
\noindent
and we know that exponentiation to $d_1$ in $K$ induces multiplication by $d_1$ in cohomology with compact support. Thus $\delta$ gives rise to the following commutative diagram with exact rows in cohomology with compact support:
$$\begin{array}{ccccccccccc}
0&&\bz/r\bz&&\bz/r\bz&&0&\\
\|&&\wr|&\simeq&\wr|&&\|&\\
\hc^1(\{0\})&\to&\hc^2(K^{\ast})&\to&\hc^2(K)&\to&\hc^2(\{0\})\\
\\
&&\downarrow\kappa_1&&\downarrow\cdot d_1&&&&\\
\\
\hc^1(\{0\})&\to&\hc^2(K^{\ast})&\to&\hc^2(K)&\to&\hc^2(\{0\})\\
\|&&\wr|&\simeq&\wr|&&\|&\\
0&&\bz/r\bz&&\bz/r\bz&&0&\\
\end{array}$$
\noindent
Since by assumption $d_1$ and $r$ are not coprime, multiplication by $d_1$ in $\bz/r\bz$ is not injective. It follows that $\kappa_1$ is not injective. Now let $n>1$ and suppose the claim true for all smaller $n$. Without loss of generality we may assume that $d_1$ is not prime to $r$. Note that $\{0\}\times (K^{\ast})^{n-1}$ is a closed subset of $K\times (K^{\ast})^{n-1}$ and $(K\times (K^{\ast})^{n-1})\setminus(\{0\}\times (K^{\ast})^{n-1})=(K^{\ast})^n$. After identifying $\{0\}\times (K^{\ast})^{n-1}$ with $(K^{\ast})^{n-1}$ we have an exact sequence of cohomology with compact support
\begin{equation}\label{exact}\hc^n(K\times(K^{\ast})^{n-1})\rightarrow\hc^n((K^{\ast})^{n-1})\rightarrow\hc^{n+1}((K^{\ast})^n)\rightarrow\hc^{n+1}(K\times (K^{\ast})^{n-1})\end{equation}
\noindent
According to the K\"unneth formula, for all indices $i$ we have
$$\hc^i(K\times(K^{\ast})^{n-1})\simeq\bigoplus_{p+q_1+\cdots+q_{n-1}=i}\hc^p(K)\otimes\hc^{q_1}(K^{\ast})\otimes\cdots\otimes\hc^{q_{n-1}}(K^{\ast}).$$
\noindent
Since $\hc^p(K)=0$ for $p\ne2$ and $\hc^q(K^{\ast})=0$ for $q\neq1,2$, it follows that 
$$\hc^n(K\times(K^{\ast})^{n-1})=0.$$
\noindent Let $$\delta':(K^{\ast})^{n-1}\rightarrow  (K^{\ast})^{n-1}$$
$$(u_1,\dots, u_{n-1})\mapsto (u_1^{d_1},\dots, u_{n-1}^{d_{n-1}})$$
\noindent be the map obtained from $\delta$ by restriction,   and 
let  $$\kappa_{n-1}:\hc^{n}((K^{\ast})^{n-1})\to\hc^{n}((K^{\ast})^{n-1}).$$
\noindent be the homomorphism it induces in cohomology with compact support. Then (\ref{exact}) gives rise to the following commutative diagram with exact rows:
$$\begin{array}{ccccc}
0&\rightarrow&\hc^n((K^{\ast})^{n-1})&\rightarrow&\hc^{n+1}((K^{\ast})^n)\\
&&\downarrow\kappa_{n-1}&&\downarrow\kappa_n\\
0&\rightarrow&\hc^n((K^{\ast})^{n-1})&\rightarrow&\hc^{n+1}((K^{\ast})^n)\\
\end{array}
$$
\noindent
Since, by induction, $\kappa_{n-1}$ is not injective, it follows that $\kappa_n$ is not injective. This completes the proof.
\par\bigskip\noindent
The following result is due to Newstead and is quoted from \cite{BS}.
\begin{lemma}\label{Newstead}Let
$W\subset\tilde W$ be affine varieties. Let $d=\dim\tilde
W\setminus W$. If there are $s$ equations $F_1,\dots, F_s$ such
that $W=\tilde W\cap V(F_1,\dots,F_s)$, then 
$$\het^{d+i}(\tilde W\setminus W,{\bz}/r{\bz})=0\quad\mbox{ for all
}i\geq s$$ and for all $r\in{\bz}$ which are prime to {\rm char}\,$K$.
\end{lemma}
From this we derive the following cohomological criterion on the number of equations defining an affine variety set-theoretically:
\begin{corollary}\label{corollary} Let $V$ be a subvariety of the affine space $K^N$. Let $s<2N$ be a positive integer such that
$$\hc^s(V,\bz/r\bz)\neq0$$
\noindent for some positive integer $r$ prime to char\,$K$. Then ara\,$V\geq N-s$. 
\end{corollary} 
\demo As in the preceding proofs, we shall omit the coefficient group $\bz/r\bz$. We have an exact sequence
$$\hc^s(K^N)\rightarrow\hc^s(V)\rightarrow\hc^{s+1}(K^N\setminus V),$$
\noindent
where $\hc^s(K^N)=0$. If $\hc^s(V)\neq0$, it follows that $\hc^{s+1}(K^N\setminus V)\ne0$. Applying Poincar\'e Duality yields 
$$\het^{2N-s-1}(K^N\setminus V)\simeq\hc^{s+1}(K^N\setminus V)\ne0,$$
\noindent
so that, by Lemma \ref{Newstead}, $V$ is not defined set-theoretically by $N-s-1$ equations, i.e., ara\,$V\geq N-s$, as claimed. This completes the proof.
\section{The cohomological lower bound}
We are now ready to complete the proof of Theorem \ref{main} for the toric variety $V$ introduced above. We show 
\begin{proposition}\label{ara} ara\,$V\geq n+1$.
\end{proposition}
\demo Let $r\in\{p,q\}$.  Suppose that char\,$K\ne r$.   It suffices to show that the claim holds under this assumption: since $p\ne q$, the characteristic of any ground field is different from $p$ or different from $q$. All the cohomology groups considered in this proof are referred to the coefficient group $\bz/r\bz$.  
Let 
$$\phi:K^n\rightarrow V$$
$$(u_1,\dots, u_n)\mapsto (u_1^{d_1}, \dots, u_{n-1}^{d_{n-1}}, u_n, u_1^{f_1}u_n^{g_1}, \dots, u_{n-1}^{f_{n-1}}u_n^{g_{n-1}}, u_1^{h_1}\cdots u_{n-1}^{h_{n-1}})$$
\noindent
and let $W\subset K^n$ be the subvariety defined by $u_n=u_1\cdots u_{n-1}=0$. The restriction
$$\phi':K^n\setminus W\rightarrow V\setminus\phi(W)$$
\noindent
is a bijection. Surjectivity is obvious, we only have to prove injectivity. Let ${\bf u}=(\bar u_1, \dots, \bar u_{n-1}, \bar u_n)\in K^n\setminus W$. First assume that $\bar u_n\ne0$. Let $i\in\{1,\dots, n-1\}$. We prove that $\bar u_i$ is uniquely determined by $\phi({\bf u})$. This is evidently true if $\bar u_i=0$. So suppose that $\bar u_i\ne0$. By (\ref{1}) there are $a_i, b_i\in{\bf Z}$ such that $a_id_i+b_if_i=1$.  Hence
$$\bar u_i=\bar u_i^{a_id_i+b_if_i}=(\bar u_i^{d_i})^{a_i}\frac{(\bar u_i^{f_i}\bar u_n^{g_i})^{b_i}}{\bar u_n^{g_ib_i}}.$$
\noindent
Now assume that $\bar u_n=0$, so that $\bar u_1,\dots, \bar u_{n-1}\ne0$. Let ${\bf v}=(\bar v_1,\dots, \bar v_{n-1},0)\in K^n$ be such that $\phi'({\bf v})=\phi'({\bf u})$. Then, for  all $i=1,\dots, n-1$, $\bar v_i^{d_i}=\bar u_i^{d_i},$
\noindent so that
\begin{equation}\label{vi}\bar v_i=\eta_i\bar u_i\end{equation}
\noindent
for some $d_i$-th root of unity $\eta_i$. We also have that
\begin{equation}\label{uv} \bar v_1^{h_1}\cdots\bar v_{n-1}^{h_{n-1}}=\bar u_1^{h_1}\cdots\bar u_{n-1}^{h_{n-1}}\end{equation}
\noindent
Since $\bar u_i\ne0$ for all $i=1,\dots, n-1$, from (\ref{vi}) and (\ref{uv}) it follows that
$$\eta_1^{h_1}\cdots\eta_{n-1}^{h_{n-1}}=1.$$
\noindent
Hence
\begin{equation}\label{omega}\omega_1=\eta_1^{h_1}=\eta_2^{-h_2}\cdots\eta_{n-1}^{-h_{n-1}}\end{equation}
\noindent
is a both a $d_1$-th root and an $\ell_1$-th root of unity, where $\ell_1=\,$lcm\,$(d_2,\dots, d_{n-1})$.  Since, by (\ref{2}), $\gcd(\ell_1,d_1)=1$, it follows that $\omega_1=1$. Hence, in view of (\ref{omega}), $\eta_1$ is both a $d_1$-th and a $h_1$-th root of unity. By (\ref{1b}) it follows that $\eta_1=1$. Similarly, one shows that $\eta_i=1$ for all $i=2,\dots, n-1$. Hence, by (\ref{vi}), ${\bf v}={\bf u}$, which completes the proof that $\phi'$ is injective. Since it is finite, it is proper. Therefore, according to \cite{CK}, Lemma 3.1, it induces an isomorphism in cohomology
$$\phi'_i:\hc^i(V\setminus\phi(W))\simeq\hc^i(K^n\setminus W),$$
for all indices $i$. Now consider the restriction map
$$\phi'':W\rightarrow \phi(W)$$
$$(u_1,\dots, u_{n-1},0)\mapsto(u_1^{d_1},\dots, u_{n-1}^{d_{n-1}},0,\dots, 0)$$
\noindent and the maps it induces in cohomology
$$\phi''_i:\hc^i(\phi(W))\rightarrow\hc^i(W).$$
\noindent
We have the following commutative diagram with exact rows:
$$\begin{array}{ccccccccccc}
&&&\alpha&&&&&\\
\hc^{n-1}(V)&\to&\hc^{n-1}(\phi(W))&\to&\hc^n(V\setminus\phi(W))&&\\
\\
&&\downarrow\phi''_{n-1}&&\wr|\downarrow\phi'_n&&&&\\
\\
\hc^{n-1}(K^n)&\to&\hc^{n-1}(W)&\to&\hc^n(K^n\setminus W)&\to&\hc^{n}(K^n)\\
\|&&&\simeq&&&\|&\\
0&&&&&&0&\\
\end{array}$$
\noindent
Our next aim is to show that $\phi''_{n-1}$ is not injective. It will follow that $\phi'_n\alpha$ is not injective, and, consequently, $\alpha$ is not injective, so that 
$$\hc^{n-1}(V)\ne0.$$
\noindent
In view of Corollary \ref{corollary}, this will imply the claim. Let $X$ be the subvariety of $K^{n-1}$ defined by $u_1\cdots u_{n-1}=0$. Then, up to omitting zero coordinates, $\phi''$ can be viewed as the map:
$$\phi'':X\rightarrow X$$
$$(u_1,\dots, u_{n-1})\mapsto (u_1^{d_1},\dots, u_{n-1}^{d_{n-1}}).$$
\noindent We have to show that the induced map 
$$\phi''_{n-1}:\hc^{n-1}(X)\rightarrow\hc^{n-1}(X)$$
is not injective. Note that $K^{n-1}\setminus X=(K^{\ast})^{n-1}$. Thus we have the following commutative diagram with exact rows, where $\kappa_{n-1}$ is the map defined in Lemma \ref{lemma1}:
$$\begin{array}{cccccccccc}
0&&&&&&0\\
\|&&&\simeq&&&\|\\
\hc^{n-1}(K^{n-1})&\to&\hc^{n-1}(X)&\to&\hc^n((K^{\ast})^{n-1})&\to&\hc^n(K^{n-1})\\
\\
&&\downarrow\phi''_{n-1}&&\downarrow\kappa_{n-1}&&&\\
\\
\hc^{n-1}(K^{n-1})&\to&\hc^{n-1}(X)&\to&\hc^n((K^{\ast})^{n-1})&\to&\hc^n(K^{n-1})\\
\|&&&\simeq&&&\|\\
0&&&&&&0\\
\end{array}$$
\noindent  By Lemma \ref{lemma1}, since $r$  is a prime factor of some of the exponents $d_i$,  if char\,$K\ne r$, $\kappa_{n-1}$ is not injective, hence nor is $\phi''_{n-1}$, as  was to be shown.  This completes the proof.
\section{Final remarks}
For the toric varieties $V$ considered in this paper, the minimal number of defining equations coincides with the minimal number of {\it binomial} defining equations, i.e., ara\,$V=$\,bar\,$V=n+1$, where ``bar'' denotes the so-called {\it binomial arithmetical rank}, a notion introduced by Thoma \cite{T2}. In general, ara\,$V<$\,bar\,$V$: an example is given by the projective monomial curves studied by Robbiano and Valla \cite{RV}, which are defined by two equations (and thus are set-theoretic complete intersections), but one of these equations is, in general, non-binomial. As it was shown in \cite{BMT2}, over fields of characteristic zero we have ara\,$V=$\,bar\,$V=\codim V$ for a toric variety $V$  only in a very special case, namely if and only if the defining ideal $I(V)$ is a complete intersection (i.e., if and only if it is generated by $\codim\,V$ equations.) The criterion changes completely in characteristic $p$: there the above equality holds for a much larger class of toric varieties, whose attached affine semigroup $\bn T$ is {\it completely p-glued}, a combinatorial property described in \cite{BMT2}. This difference explains why there are several examples of toric varieties whose arithmetical rank strictly depends on the characteristic of the ground field. One can find classes for which ara\,$V=$\,bar\,$V=\codim V$ in one positive characteristic, whereas ara\,$V>\codim V$ in all remaining characteristics: this is the case for the simplicial toric varieties of codimension 2 described in \cite{BL}, and for certain Veronese varieties, of arbitrarily high codimension, considered in \cite{B2}. In this respect, the main result of this paper presents a new situation: the arithmetical rank is constant in all characteristics. Moreover the defining equations are the same over all fields. This is not true for monomial curves in the three-dimensional projective space: in every positive characteristic they are binomial set-theoretic complete intersections, but the two defining binomial equations change from one characteristic to the other (see \cite{Mo}); whereas in characteristic zero their arithmetical ranks are unknown in general.

\end{document}